\documentclass[11pt,a4paper]{article}
\usepackage{setspace}

\usepackage[T1]{fontenc}
\usepackage[utf8]{inputenc}
\usepackage{authblk}
\usepackage[top=20mm,bottom=20mm, left=15mm,right=15mm]{geometry}
\usepackage{amsmath,amsfonts,amsthm,amssymb, graphicx,color,bm,mathrsfs,psfrag}
\usepackage[unicode,breaklinks=true,colorlinks=true]{hyperref}
\usepackage[dvipsnames]{xcolor}

\newtheorem{theorem}{Theorem}[section]
\newtheorem{lemma}[theorem]{Lemma}

\theoremstyle{definition}
\newtheorem{definition}[theorem]{Definition}

\theoremstyle{remark}
\newtheorem{remark}[theorem]{Remark}

\newtheorem*{claim*}{Claim}

\numberwithin{equation}{section}

\allowdisplaybreaks

\title{$p$-Kirchhoff type equation with Neumann boundary conditions\thanks{This work was supported by  National Natural Science Foundation of China (Grant No. 11771423, 11871452) and National Science Foundation of Jiangsu Higher Education Institutions of China (Grant No. 19KJD100007).}
}
\author[a,b]{Weihua Wang\thanks{Corresponding author: wangvh@163.com, wangweihua15@mails.ucas.ac.cn, ORCID iD: https://orcid.org/0000-0002-8802-743X
}}
\affil[a]{School of Mathematical Sciences, Yangzhou University, Yangzhou, 225002,  China}
\affil[b]{University of Chinese Academy of Sciences, Beijing, 100049,  China}

\date{} 

\begin{document}
  \maketitle
\linespread{1.0}
\begin{abstract}
This paper is concerned with the multiplicity results  to a class of  $p$-Kirchhoff type elliptic equation with the homogeneous  Neumann boundary conditions by an abstract linking lemma due to Br\'{e}zis and Nirenberg. We obtain the twofold results in subcritical and critical cases,  which is a meaningful addition and completeness to the known results about Kirchhoff  equation.

\end{abstract}
\section{Introduction}
    In this paper, we inspect the multiplicity results  for  $p$-Kirchhoff type elliptic equation with the homogeneous Neumann boundary conditions
\begin{equation}\label{KPN}\tag{KN}
\left\{\begin{array}{ll}
-\left(a+b\int_{\Omega}|\nabla u|^{p} \mathrm{d} x\right)^{s\cdot{\rm sgn}b} \triangle_{p} u = f(x,u) & \text { in } \Omega, \\
\frac{\partial u}{\partial n}= 0 & \text { on } \partial \Omega,
\end{array}\right.
 \end{equation}
where $\Omega$  is a bounded smooth domain in $\mathbb{R}^N$, $\triangle_{p}$ denotes $p$-Laplace operator,  $\frac{\partial}{\partial n}$ is the outer unit normal derivative to the boundary $\partial \Omega$,  $a, b \geq 0$,  $a+b>0$  and $0<p(s\cdot{\rm sgn}b + 1) \leq p^{*} := \frac{Np}{N-p}$.\\
When $b=0$, Eq. \eqref{KPN} is reduced to a $p$-Laplacian that is extensively studied and widely used.\\
When $s=1, \, p=2,\, \&\,  b>0$, Eq. \eqref{KPN} is  the stationary version of  the second order hyperbolic equation
\begin{equation}\label{Kir}\tag{Kir}
 \rho \frac{\partial^{2} u}{\partial t^{2}}-\left(\frac{P_{0}}{h}+\frac{E}{2 L} \int_{0}^{L}\left|\frac{\partial u}{\partial x}\right|^{2} d x\right) \frac{\partial^{2} u}{\partial x^{2}}=0
\end{equation}
presented by Kirchhoff in \cite{Ki76} to describe  a string vibration equation that does not ignore the change in string's length during the vibration. A obvious characteristic is that the \eqref{Kir} involve a nonlocal coefficient $\frac{P_{0}}{h}+\frac{E}{2 L} \int_{0}^{L}\left|\frac{\partial u}{\partial x}\right|^{2} d x$ which rely on the average $\frac{1}{2 L} \int_{0}^{L}\left|\frac{\partial u}{\partial x}\right|^{2} d x,$ and therefore the equation is no more a point wise identity. The parameters in \eqref{Kir} have the following meanings: $L$ denotes the length of the string, $h$ represents the area of the cross-section, $E$ means the Young modulus of the material, $\rho$ and $P_0$ is the mass density and  the initial tension, respectively.
And this type of problem also comes from biological systems that $u$ represents a process that depends on its own averaging\cite{ACM05,CR92} and references therein.

$p$-Kirchhoff type elliptic equations such as Eq. \eqref{KPN}  have been extensively studied by the  nonlinear functional analysis approach since Lions' work\cite{Lions78}. The main goal of these works are to study the effect of the nonlocal coefficients: $M(\int_{\Omega}|\nabla u|^{p} d x)$ on the principal term of these equations compared with the common second-order elliptic equations. Naimen \cite{Na14} study the existence and nonexistence of solutions of $-\left(a+b \int_{\Omega}|\nabla u|^{2} d x\right) \Delta u=\mu g(x, u)+u^{5}, u>0$ in $H^1_{0}(\Omega)$, which is some extension of a part of Brezis-Nirenberg's result\cite{BN83}. J\'{u}lio etc. \cite{JCF06} considered  $-\left[M\left(\int_{\Omega}|\nabla u|^{p} d x\right)\right]^{p-1} \triangle_{p} u=f(x, u)+\lambda|u|^{s-2} u$ in $W_{0}^{1,p}(\Omega)$ with crtical/supercritical growth by using variational methods combined with Moser's iteration method under the condition that $M(t)\geq m_0 >0$ and $M(t)= m_1 \geq m_0$ for all $t\geq t_0$ or its equivalent variational condition. The above restrictions are removed from the results in our article.  More about the existence and multiplicity of solutions for Kirchhoff/$p$-Kirchhoff (type) equations, we refer to \cite{AA17,Ch18,FOS19,MZ09,PK06,Wang19} and therein.

 Chabrowski\cite{Ch18} studied  the existence of positive solutions for
$$-\left(\int_{\Omega}|\nabla u|^{2} d x\right)^{s} \triangle u=Q(x)|u|^{\ell-2} u+\left(\int_{\Omega}|u|^{q} d x\right)^{r}|u|^{q-2} u$$ with Neumann boundary conditions. Motivated by \cite{Ch18}, we deal with $p$-Kirchhoff type equation and the classical $p$-Laplacian in a uniform form in this article. We will show that $p$-Kirchhoff type equation can still be treated like the semilinear elliptic equations with $p$-Laplacian despite the presence of the nonlocal coefficients. In addition,  we obtain the twofold results in not only subcritical but also critical cases.

 Boureanu and Mih\u{a}ilescu \cite{BM08} established the existence and multiplicity of solutions for $p(x)$-Laplacian with Neumann boundary conditions in case of subcritical growth, which  is the main source of our method in the present article.

 For the purpose of  depicting our results, we need to present the concept of weak solution for Eq. \eqref{KPN} and the right side term $f$ satisfied conditions:
 \begin{definition}
We say that $u \in W^{1,p}(\Omega)$ is a weak solution of  Eq. \eqref{KPN} provided
$$\left(a +b\int_{\Omega}|\nabla u|^p dx\right)^{s\cdot{\rm sgn}b}\int_{\Omega}|\nabla u|^{p-2}\nabla u \cdot \nabla\varphi dx - \int_{\Omega}f(x,u)\varphi dx =0$$
for any $\varphi$ in $W^{1,p}(\Omega)$.

\end{definition}

The nonlinear term $f : \overline{\Omega} \times \mathbb{R} \rightarrow \mathbb{R}$ is assumed to be a continuous function that satisfies the following  condition:
\begin{description}
  \item[(f$_1$)] (subcritical growth)\, $|f(x, t)| \leq C\left(1+|t|^{q-1}\right)$ for all $t \in \mathbb{R}$ and $x \in \overline{\Omega}$, where $C$ is a positive constant and $p(s\cdot{\rm sgn}b + 1) < q < p^{*}$.
  \item[(\~{f}$_1$)] (critical growth)\, $f(x, t)= -|t|^{p^{*}-2}t + g(x,t)$, where  $g(x, t)$ satisfies the following conditions:\\
      $|g(x, t)| \leq \tilde{C}|t|^{q-1}$ for all $t \in \mathbb{R}$ and $x \in \overline{\Omega}$ with $\tilde{C}$ is a positive constant and $p(s\cdot{\rm sgn}b + 1) \leq q < p^{*}$.
  \item[(f$_2$)] There exists a constant $\eta_{0}>0$ such that $0<\lim _{|t| \rightarrow 0} \frac{F(x, t)}{|t|^{p(s\cdot{\rm sgn}b+1)}}<\eta_{0}$ uniformly a.e. $x \in \Omega$.
  \item[(f$_3$)] $\lim _{|t| \rightarrow+\infty} F(x, t)=-\infty$ and $\lim _{|t| \rightarrow \infty} \frac{F(x, t)}{|t|^{p(s\cdot{\rm sgn}b+1)}} \leq 0$ uniformly a.e.$x \in \Omega$, where $F(x, t)=\int_{0}^{t} f(x, s) d s$.
\end{description}
\begin{remark}
The hypothesis (\~{f}$_1$) can be weakened to (\b{f}$_1$)\\ $f(x, t)= -|t|^{p^{*}-2}t + g(x,t)$ with  $g(x, t)$ satisfying the conditions:\\
      \begin{itemize}
        \item $|g(x, t)| \leq \tilde{C}(1+|t|^{q-1})$  for all $t \in \mathbb{R}$ and $x \in \overline{\Omega}$,
        \item $\lim _{|t| \rightarrow \infty} \frac{g(x, t)}{|t|^{p(s\cdot{\rm sgn}b+1)-1}} \leq 0$ uniformly a.e.$x \in \Omega$,
      \end{itemize}
         where $\tilde{C}$ is a positive constant and $p(s\cdot{\rm sgn}b + 1) \leq q < p^{*}$.
\end{remark}

For the case of  $f(x, t)$ with subcritical growth, we have the following results.
\begin{theorem}[Subcritical case]\label{Thm1}
Assume that the function $f(x,t)$ satisfies (f$_1$), (f$_2$) and (f$_3$). Then problem \eqref{KPN} owns at least two non-trivial weak solutions.
\end{theorem}

For the case of  $f(x, t)$ with critical growth, we have the following results.
\begin{theorem}[Critical case ]\label{Thm2}
Suppose that the function $f(x,t)$ satisfies (\~{f}$_1$). Then problem \eqref{KPN} possesses at least two non-trivial weak solutions if $a>0$.
\end{theorem}

This paper is organized as follows: In Section 2, we review some necessary preliminaries. In Section 3, we gets the results
when $f(x, t)$ is subcritical growth on $t$. And in Section 4, we deal with the case of critical growth.

\section{Preliminaries}
We define the energy functional $I_{b}(u): W^{1,p}(\Omega)\rightarrow\mathbb{R}$ as
\begin{equation*}
  I_{b}(u)=\left\{\begin{array}{ll}
  \frac{1}{p(s+1)b}\left(a +b\int_{\Omega}|\nabla u|^p dx\right)^{s+1}  - \frac{a^{s+1}}{p(s+1)b} - \int_{\Omega}F(x,u) dx, & b>0,\\
  \frac{1}{p}\int_{\Omega}|\nabla u|^p dx - \int_{\Omega}F(x,u) dx, & b=0.
  \end{array}\right.
\end{equation*}
Standard arguments indicate that $I_{b}\in C^1(W^{1,p}(\Omega), \mathbb{R})$ with
\begin{equation*}
  I_{b}'(u)\varphi=\left(a +b\int_{\Omega}|\nabla u|^p dx\right)^{s\cdot{\rm sgn}b}\int_{\Omega}|\nabla u|^{p-2}\nabla u \cdot\nabla\varphi dx - \int_{\Omega}f(x,u)\varphi dx,
\end{equation*}
for any $u, \varphi$ in $W^{1,p}(\Omega)$. Therefore, weak solutions of Eq. \eqref{KPN} are just the critical points of $I_{b}$.

We  also need  Palais-Smale  "compactness" condition:
\begin{definition}[Palais, 1970, Definition 1.3, \cite{GrT01}]
Let $X$ be a Banach space. A $C^{1}$-functional $J : X \rightarrow \mathbb{R}$ satisfies Palais-Smale (henceforth denoted by (PS)) condition if every sequence $\left\{u_{n}\right\}$ in $X$ such that $J\left(u_{n}\right)$ is bounded and $\lim _{n\to +\infty} J^{\prime}\left(u_{n}\right)=0$ in $X^{*}$ has a convergent subsequence.
\end{definition}
 A basic tool in this paper is the following abstract linking arguments due to Br\'{e}zis and Nirenberg\cite{BN91}.
\begin{lemma}[\cite{BN91}, Theorem 4]\label{Lem2.1}
Suppose that $\left(X,\|\cdot\|_{X}\right)$ is a Banach space with the direct sum decomposition $X=X_{1} \oplus X_{2}$ with $\operatorname{dim}\left(X_{2}\right)<+\infty .$ Assume that $J \in C^{1}(X, \mathbb{R})$ with $J(0)=0$ satisfies (PS) condition in $X$. Moreover, for a  constant $\rho>0,$  we have
$$J(u) \geq 0 \text { for all } u \in X_{1} \text { with }\|u\|_{X} \leq \rho$$ and $$J(u) \leq 0 \quad \text { for all } u \in X_{2} \text { with }\|u\|_{X} \leq \rho.$$
Also assume that $J$ is bounded from below and inf$_{X} J < 0$. Then $J$ has at least two nontrivial critical points.
\end{lemma}

In order to apply Lemma \ref{Lem2.1} to the functional $I_{b}$, we have to decompose the space $W^{1,p}(\Omega)$ as $W^{1,p}(\Omega)=W_{0} \oplus \mathbb{R}$, where  $$W_{0}=\left\{u \in W^{1,p}(\Omega) : \int_{\Omega} u d x=0\right\}.$$
Obviously, $W_{0}$ is a closed linear subspace of $W^{1,p}(\Omega)$ with codimension $1$. In detail, under the above space decomposition, for any $u\in W^{1,p}(\Omega)$, denote $\overline{u}=\frac{1}{|\Omega|} \int_{\Omega} u d x $ and $\widetilde{u}=u-\overline{u}$. Then $u=\widetilde{u}+\overline{u}$,  where $\overline{u} \in \mathbb{R}$ and $\widetilde{u} \in W_{0}$. This decomposition yields the norm $\|u\|^{p}=\int_{\Omega}|\nabla \widetilde{u}|^{p} d x + \overline{u}^{p}$ of $W^{1,p}(\Omega)$, which is equivalent to the "standard" norm $\|u\|_{W^{1,p}(\Omega)}:=\|\nabla u\|_{L^p(\Omega)}+\|u\|_{L^p(\Omega)}$ on $W^{1,p}(\Omega)$. For the sake of simplicity, we will not distinguish between these two norms in the following discussion.

In order to get the functional $I_{b}$ to meet the (PS) condition, we also need the following inequality:
\begin{lemma}[Lemma 5.1 and Lemma 5.2, \cite{GM75} or Lemma 3.3, \cite{JCF06}]\label{Lem2.2}
Let  $\langle \cdot , \cdot\rangle$ be the normal inner product in $\mathbb{R}^{N}.$ Then
$$\left\langle|x|^{p-2} x-|y|^{p-2} y, x-y\right\rangle \geq C_{p}|x-y|^{p} \text { if } p \geqslant 2$$
or
$$\left\langle|x|^{-2} x-|y|^{p-2} y, x-y\right\rangle \geq \frac{C_{p}|x-y|^{2}}{(|x|+|y|)^{2-p}} \text { if } 2>p>1$$  for any $x, y$ be in $\mathbb{R}^{N}$.
\end{lemma}

\section{Subcritical case}

\begin{lemma}\label{Lem3.1}
Assume that conditions (f$_1$) and (f$_3$) are satisfied.  Then the functional $I_{b}$ fulfills the $(P S)$ condition.
\end{lemma}
\begin{proof}
Let $\left\{u_{n}\right\} \subseteq W^{1,p}(\Omega)$ be such that $|I_{b}(u_{n})| \leq M$ and $\left|\langle I_{b}^{\prime}(u_{n}), \varphi\rangle\right| = o(\|\varphi\|)$ for any $ \varphi \in W^{1,p}(\Omega)$.

 {\bf\it Claim}:  $\left\{u_{n}\right\}$ is bounded in $W^{1,p}(\Omega)$. \\
 In fact, arguing by contradiction and passing to a subsequence, we assume that $\left\|u_{n}\right\| \rightarrow \infty$ as $n\to +\infty$.

Let $v_{n}(x) :=\frac{u_{n}(x)}{\left\|u_{n}\right\|}$, then $\|v_{n}\|=1$ in $W^{1,p}(\Omega)$ and $W^{1,p}(\Omega)$ is a reflexive Banach space, up to a
subsequence, $v_n \rightharpoonup v$ in $W^{1,p}(\Omega)$. Therefore,
\begin{eqnarray*}
 v_n&\rightarrow& v  \text{ in } L^r(\Omega), 1< r< p^*,\\
v_n(x)&\rightarrow& v(x) \,\, a.e. \text{ in } \Omega.
\end{eqnarray*}
 Since $\lim _{|t| \rightarrow \infty} \frac{F(x, t)}{|t|^{p(s\cdot{\rm sgn}b+1)}} \leq 0$, for any $\varepsilon >0$, there exists $R >0$ such that $F(x,t)\leq \varepsilon |t|^{p(s\cdot{\rm sgn}b+1)}$ for all $(x, |t|)$ in $\Omega\times [R, +\infty)$. And together with $F$ is continuous in $\overline{\Omega} \times \mathbb{R}$, there exists a constant $C_1 >0$ such that
 \begin{equation}\label{Eq3.1}
   F(x,t)\leq \varepsilon |t|^{p(s\cdot{\rm sgn}b+1)} +C_1
 \end{equation}
 for any $(x, |t|)$ in $\Omega\times [R, +\infty)$.\\

 First of all, we  consider the case of $b>0$.
 \begin{eqnarray}\label{Eq3.2}
   \frac{M}{\|u_n\|^{p(s+1)}} &\geq& \frac{I_{b}(u_n)}{\|u_n\|^{p(s+1)}} \nonumber\\
                                &=&  \frac{1}{\|u_n\|^{p(s+1)}}\left[ \frac{1}{p(s+1)b}\left(a +b\int_{\Omega}|\nabla u_n|^p dx\right)^{s+1}  - \frac{a^{s+1}}{p(s+1)b} - \int_{\Omega}F(x,u_n) dx \right]\nonumber\\
                                &\geq&  \frac{1}{p(s+1)b}\left(\frac{a}{\|u_n\|^{p}} +b\int_{\Omega}|\nabla v_n|^p dx\right)^{s+1}  - \frac{a^{s+1}}{p(s+1)b\|u_n\|^{p(s+1)}} - \varepsilon C_2- \frac{C_1|\Omega|}{\|u_n\|^{p(s+1)}},
 \end{eqnarray}
where $C_2$ is positive constant.

By \eqref{Eq3.2} and $\lim_{n\to +\infty}\|u_n\|= +\infty$, together with the arbitrarily of $\varepsilon$, we have $\lim_{n\to +\infty}\int_{\Omega}|\nabla v_n|^p dx = 0$.
Therefore, $$0\leq \int_{\Omega}|\nabla v_0|^p dx \leq \liminf_{n\to +\infty}\int_{\Omega}|\nabla v_n|^p dx = 0.$$
Hence, $\nabla v_0 = 0$ a.e. $x\in \Omega$ which yields $v_0$ in $\mathbb{R}$ and
\begin{equation}\label{Eq3.3}
  \lim_{n\to +\infty}\int_{\Omega}|\nabla (v_n - v_0)|^p dx =\lim_{n\to +\infty}\int_{\Omega}|\nabla v_n|^p dx =0.
\end{equation}
Considering $v_n \rightarrow v  \text{ in } L^r(\Omega), 1< r< p^*$, we have $v_n\rightarrow v  \text{ in } W^{1,p}(\Omega)$. That fact combined
with $\left\|v_{n}\right\|=1$ shows that $v_{0} \neq 0$ and consequently
$$\lim_{n\to \infty}|u_{n}(x)|= \lim_{n\to \infty}\|u_{n}(x)\|\cdot|v_{n}|=\lim_{n\to \infty}\|u_{n}(x)\|\lim_{n\to \infty}|v_{n}| = +\infty, a.e. x \in \Omega.$$  Hence, $\lim _{n \rightarrow \infty} \int_{\Omega} F\left(x, u_{n}\right) d x=-\infty$ by (f$_3$).
Therefore,
\begin{eqnarray*}
  M &\geq& I_{b}(u_{n})= \frac{1}{p(s+1)b}\left(a +b\int_{\Omega}|\nabla u_{n}|^p dx\right)^{s+1}  - \frac{a^{s+1}}{p(s+1)b} - \int_{\Omega}F(x,u_{n}) dx\\
     &\geq&  - \frac{a^{s+1}}{p(s+1)b} - \int_{\Omega}F(x,u_{n}) dx \to +\infty \text{ as } n\to +\infty,
\end{eqnarray*}
which is absurd. Hence, $u_{n}$ is bounded in $W^{1,p}(\Omega)$ in case of $b>0$.\\
 And we can use a similar method to prove that $u_n$ is bounded in $W^{1,p}(\Omega)$ when $b=0$. Therefore, $u_{n}$ is bounded in $W^{1,p}(\Omega)$.
 Passing to a subsequence if necessary, there exists $u \in W^{1,p}(\Omega)$ such that  $u_n \rightharpoonup u$ in $W^{1,p}(\Omega)$ and thus
\begin{eqnarray}\label{Eq3.4}
 u_n&\rightarrow& u  \text{ in } L^r(\Omega) \text{ with } 1< r< p^*,\nonumber\\
u_n(x)&\rightarrow& u(x) \,\, a.e. \text{ in } \Omega.
\end{eqnarray}
 According to \eqref{Eq3.4} and the condition (f$_1$),  we have
\begin{eqnarray}\label{Eq3.5}
   & & \left|\int_{\Omega} f\left(x, u_{n}\right)\left(u_{n}-u\right) d x\right| \nonumber\\
   &\leq& C \int_{\Omega}\left(1+\left|u_{n}\right|^{q-1}\right)\left|u_{n}-u\right| d x \nonumber\\
   &\leq& C_{1}\|u_{n}\|^{q-1}_{q^{\prime}}\left\|u_{n}-u\right\|_{q}+C_{1}|\Omega|^{\frac{1}{p^{\prime}}}\left\|u_{n}-u\right\|_{p},
\end{eqnarray}with $q^{\prime}=\frac{q}{q-1}$ and $p^{\prime}=\frac{p}{p-1}$,
which approaches 0 as $n \rightarrow \infty$.

In the light of $\left|\langle I_{b}^{\prime}(u_{n}), \varphi\rangle\right| = o(\|\varphi\|)$ and \eqref{Eq3.5}, we have
\begin{equation*}
  I_{b}'(u_{n})(u_{n}-u)=\left(a +b\int_{\Omega}|\nabla u_{n}|^p dx\right)^{s\cdot{\rm sgn}b}\int_{\Omega}|\nabla u_{n}|^{p-2}\nabla u_{n} \cdot\nabla(u_{n}-u) dx - \int_{\Omega}f(x,u_{n})(u_{n}-u) dx.
\end{equation*}
Hence,
\begin{eqnarray*}
   & & \lim_{n\to +\infty} \left(a +b\int_{\Omega}|\nabla u_{n}|^p dx\right)^{s\cdot{\rm sgn}b}\int_{\Omega}|\nabla u_{n}|^{p-2}\nabla u_{n}\cdot \nabla(u_{n}-u) dx\\
   &=& \lim_{n\to +\infty}I_{b}'(u_{n})(u_{n}-u)  + \lim_{n\to +\infty}\int_{\Omega}f(x,u_{n})(u_{n}-u) dx \\
   &=& 0+0=0.
\end{eqnarray*}
Since $u_{n}$ is bounded in $W^{1,p}(\Omega)$, hence,  $\lim_{n\to +\infty} \int_{\Omega}|\nabla u_{n}|^{p-2}\nabla u_{n} \cdot\nabla(u_{n}-u) dx =0$.

On the other hand, $u_n \rightharpoonup u$ in $W^{1,p}(\Omega)$, we infer that $$\lim_{n\to +\infty} \int_{\Omega}|\nabla u|^{p-2}\nabla u \cdot\nabla(u_{n}-u) dx =0.$$ And utilizing Lemma \ref{Lem2.2}, we may obtain
\begin{equation}\label{Eq3.6}
  \lim_{n\to +\infty} \int_{\Omega}|\nabla(u_{n}-u)|^pdx =0.
\end{equation}
Combining the above fact with \eqref{Eq3.6} and \eqref{Eq3.4}, we deduce that  $u_{n} \rightarrow u$ in $W^{1,p}(\Omega)$, i.e. $I_{b}$ appeases the $(P S)$ condition.
\end{proof}

\begin{lemma}\label{Lem3.2}
Assume that condition (f$_2$) and (f$_3$) are contented. Then $I_{b}$ is bounded from below and $\underset{u\in H^1(\Omega)}{\inf}  I_{b} (u)<0$.
\end{lemma}
\begin{proof}
Since $\lim _{|t| \rightarrow+\infty} F(x, t)=-\infty$ and $F$ is continuous in $\overline{\Omega} \times \mathbb{R}$ , we  infer that there exists a constant $M>0$ such that
$\int_{\Omega}F(x,u)dx < M_3.$
Therefore, we have
\begin{equation*}
  I_{b}(u)>\left\{\begin{array}{ll}
  \frac{1}{p(s+1)b}\left(a +b\int_{\Omega}|\nabla u|^p dx\right)^{s+1}  - \frac{a^{s+1}}{p(s+1)b} -  M_3, & b>0\\
  \frac{1}{p}\int_{\Omega}|\nabla u|^p dx -  M_3, & b=0
  \end{array}\right\}
  \geq
  \left\{\begin{array}{ll}
   - \frac{a^{s+1}}{p(s+1)b} -  M_3, & b>0,\\
   -  M_3, & b=0.
  \end{array}\right\} >-\infty,
\end{equation*}
that is, $I_{b}$ is bounded from below.

Since $0<\lim _{|t| \rightarrow 0} \frac{F(x, t)}{|t|^{p(s\cdot{\rm sgn}b+1)}}<\eta_{0}$, we deduce that there exists a constant $\delta >0$ such that, for any $|t|\in (0, \delta)$,
$0< \frac{F(x, t)}{|t|^{p(s\cdot{\rm sgn}b+1)}}<\eta_{0}$ or $0< F(x, t)<\eta_{0} |t|^{p(s\cdot{\rm sgn}b+1)}$.

Notice that $F(x,0)\equiv 0$, we obtain
\begin{equation}\label{Eq3.7}
 0 \leq F(x, t) \leq  \eta_{0} |t|^{p(s\cdot{\rm sgn}b + 1)} \text { for any } (x, |t|) \in \overline{\Omega} \times [0, \delta].
\end{equation}
Hence, there exists a constant function $t_0$ in $W^{1,p}(\Omega)$ with $0<t_0<\delta$ such that
\begin{equation*}
  I_{b}(t_0)= -\int_{\Omega}F(x,t_0) dx < 0.
\end{equation*}
Thence, we have $\underset{u\in W^{1,p}(\Omega)}{\inf}  I (u) \leq  I_{b}(t_0) < 0$.

\end{proof}

\begin{lemma}\label{Lem3.3}
Assume that conditions (f$_1$) and (f$_2$) are sufficed. Then there exists $\rho>0$ such that for all $u \in H_0$ with $\|u\| \leq \rho$ we have $I_{b}(u) \geq 0$ and $I_{b}(e) \leq 0$ for all $e \in \mathbb{R}$ with $|e| \leq \rho$.
\end{lemma}
\begin{proof}
We take $u \in W_0$ with $\|u\|=\rho,$ where $\rho$ is small enough and will be determined later.
Combining \eqref{Eq3.7} with the condition (f$_1$), there exists a constant $C_{\eta_{0}}$ dependent on $\eta_{0}$ in (f$_2$) such that
\begin{equation}\label{Eq3.8}
 F(x, t) \leq  \eta_{0} |t|^{p(s\cdot{\rm sgn}b + 1)} + C_{\eta_{0}} |t|^q  \text { for any } (x, t) \in \overline{\Omega} \times \mathbb{R}.
\end{equation}
Thence, by continuous embedding,
\begin{eqnarray*}
  \int_{\Omega} F(x,u)dx &\leq& \eta_{0}\int_{\Omega} |u|^{p(s\cdot{\rm sgn}b + 1)}dx +  C_{\eta_{0}}\int_{\Omega} |u|^q dx\\
                         &\leq&  C_2\eta_{0}\|u\|^{p(s\cdot{\rm sgn}b + 1)} +  \overline{C}_{\eta_{0}}\|u\|^q.
\end{eqnarray*}
By Poincar\'{e}'s inequality, there exists a positive constant $C_3$ such that $\|\nabla u\|_{L^p(\Omega)}\geq C_3\|u\|$.
If $b>0$,
\begin{eqnarray*}
  I_{b}(u) &=& \frac{1}{p(s+1)b}\left(a +b\int_{\Omega}|\nabla u|^p dx\right)^{s+1}  - \frac{a^{s+1}}{p(s+1)b} - \int_{\Omega}F(x,u) dx \\
           &\geq& \frac{b^s}{p(s+1)}\left(\int_{\Omega}|\nabla u|^p dx\right)^{s+1} -  C_2\eta_{0}\|u\|^{p(s + 1)} -  \overline{C}_{\eta_{0}}\|u\|^q \\
            &\geq& \left[ \frac{C^{p(s+1)}_{3}b^s}{p(s+1)} - \frac{3}{2} C_2\eta_{0}\right]\|u\|^{p(s + 1)} + \frac{1}{2} C_2\eta_{0}\|u\|^{p(s + 1)} -  \overline{C}_{\eta_{0}}\|u\|^q
\end{eqnarray*}
Take $\frac{C^{p(s+1)}_{3}b^s}{p(s+1)} - \frac{3}{2} C_2\eta_{0}=0$, i.e. $\eta_{0}= \frac{2C^{p(s+1)}_{3}b^s}{3pC_{2}(s+1)} >0$, then
\begin{equation}\label{}
  I_{b}(u) \geq \left[ \frac{C^{p(s+1)}_{3}b^s}{3p(s+1)} - \overline{C}_{\eta_{0}}\|u\|^{q - p(s + 1)}\right]\|u\|^{p(s + 1)}.
\end{equation}
Combined with the above inequalities and the condition $p(s\cdot{\rm sgn}b + 1)<q$, it yields that there exists $\rho >0$ small enough such that $I_{b}(u)\geq 0$ for all $u$ in $W_0$ with $\|u\|\leq \rho$.
And we can use a similar method to prove that $I_{b}(u)\geq 0$ for all $u$ in $W_0$ with $\|u\|\leq \rho$ when $b=0$.

At last, considering the constant function $t$ in $\mathbb{R}\subset W^{1,p}(\Omega)$, we have $I_{b}(u) = -\int_{\Omega}F(x,t)dt$, which yields $I_{b}(u)\leq 0$ for $t\in \mathbb{R}$ small enough together with \eqref{Eq3.7}.
\end{proof}

\begin{proof}[Proof of Theorem \ref{Thm1}]
Combining  Lemma \ref{Lem2.1} with   Lemma \ref{Lem3.1},  Lemma \ref{Lem3.2} and  Lemma \ref{Lem3.3}, we finish the proof of Theorem \ref{Thm1}.
\end{proof}

\section{Critical case}
In the critical case, by the hypothesis (\~{f}$_1$), we have
\begin{equation*}
  I_{b}(u)=\left\{\begin{array}{ll}
  \frac{1}{p(s+1)b}\left(a +b\int_{\Omega}|\nabla u|^p dx\right)^{s+1}  - \frac{a^{s+1}}{p(s+1)b}  + \int_{\Omega} \frac{1}{p^{*}}|u|^{p^{*}}dx -  \int_{\Omega}G(x,u)dx, & b>0,\\
  \frac{1}{p}\int_{\Omega}|\nabla u|^p dx  + \int_{\Omega} \frac{1}{p^{*}}|u|^{p^{*}}dx  -  \int_{\Omega}G(x,u)dx, & b=0,
  \end{array}\right.
\end{equation*}
where $G(x,t)=\int_{0}^{t}g(x,\tau)d\tau$.
Standard arguments show that $I_{b}\in C^1(W^{1,p}(\Omega), \mathbb{R})$ with
\begin{equation*}
  I_{b}'(u)\varphi=\left(a +b\int_{\Omega}|\nabla u|^2 dx\right)^{s\cdot{\rm sgn}b}\int_{\Omega}\nabla u \cdot\nabla\varphi dx +\int_{\Omega}|u|^{2^{*-2}}u\varphi dx - \int_{\Omega}g(x,u)\varphi dx,
\end{equation*}
for any $u, \varphi$ in $W^{1,p}(\Omega)$.

\begin{lemma}\label{Lem4.1}
Assume that the condition (\~{f}$_1$) is satisfied.  Then the function $I_{b}$ fulfills the $(P S)_c$ condition.
\end{lemma}
\begin{proof}
The condition (\~{f}$_1$) indicates (f$_3$).  From the fact $\left\{u_{n}\right\}$ is a Palais-Smale sequence, it follows that $\{u_{n}\}$ is bounded in $W^{1,p}(\Omega)$ (Refer to {\bf\it Claim} in the proof of Lemma \ref{Lem3.1}). Thus, there exists a subsequence, that we will denote by $u_{n}$ such that
\begin{eqnarray*}
  u_{n} &\rightharpoonup& u \text { weakly in } W^{1,p}(\Omega), \\
  u_{n} &\rightarrow& u \quad \text { strongly in } L^{r}(\Omega), 1<r<p^{*} \text { and for almost every } x \text { in } \Omega.
\end{eqnarray*}
The Concentration Compactness Lemma\cite{Lions85} indicates:
If $|\nabla u_{n}|^{p} \rightharpoonup \mu \text{ and } |u_{n}|^{p^{*}} \rightharpoonup \eta$ weakly-$^*$ in the sense of measures, where $\mu$ and $\eta$ are bounded nonnegative measures on $\mathbb{R}^N$, then there exist at most countable $x_{1}, x_{2}, \ldots, x_{K} \in \overline{\Omega}$ and $\eta_{1}, \eta_{2}, \ldots, \eta_{K}, \mu_{1}, \mu_{2}, \ldots, \mu_{K}$ nonnegative numbers such  that
\begin{eqnarray}\label{Eq4.1}
  \eta &=& |u|^{p^{*}}+\sum_{j=1}^{K} \eta_{j} \delta_{x_{j}}, \quad \eta_{j}>0 \nonumber\\
  \mu &\geq& |\nabla u|^{p}+\sum_{j=1}^{K} \mu_{j} \delta_{x_{j}}, \quad \mu_{j}>0 \nonumber\\
  \left(\eta_{j}\right)^{\frac{p}{p^{*}}} &\leq& \frac{\mu_{j}}{S}
\end{eqnarray}
where $\delta_{x}$ is the Dirac-mass of mass 1 concentrated at $x \in \overline{\Omega}$.

Let $\phi \in C^{\infty}(\mathbb{R}^{N})$ with $\phi \equiv 1$ in $B(x_{k}, \varepsilon) \text{ and } \phi \equiv 0$ in $B(x_{k}, 2 \varepsilon)^{c},|\nabla \phi| \leqslant \frac{2 C}{\varepsilon}$, where $x_{i} \in \overline{\Omega}$ belongs to the support of $\eta .$  Considering the boundedness of  the sequence $\left\{u_{n} \phi\right\}$  in $W^{1,p}(\Omega)$,  we obtain
\begin{eqnarray*}
 0 \leftarrow\left\langle I_{b}^{\prime}\left(u_{n}\right), u_{n} \phi\right\rangle&=&\left(a +b\int_{\Omega}|\nabla u_{n}|^p dx\right)^{s\cdot{\rm sgn}b}\int_{\Omega}|\nabla u_{n} |^{p-2}\nabla u_{n} \cdot\nabla(u_{n} \phi) dx \\
    & & \qquad + \int_{\Omega}\left|u_{n}\right|^{p^{*}} \phi d x-  \int_{\Omega}g(x,u_{n})u_{n} \phi dx.
\end{eqnarray*}
that is
\begin{eqnarray}\label{Eq4.1}
   & & - \int_{\Omega} \phi d \eta+  \int_{\Omega} g(x, u) u \phi d x \nonumber\\
   &=& \lim _{n \to \infty} \left(a +b\int_{\Omega}|\nabla u_{n}|^p dx\right)^{s\cdot{\rm sgn}b}\int_{\Omega}|\nabla u_{n} |^{p-2}\nabla u_{n} \cdot\nabla(u_{n} \phi) dx\nonumber\\
   &=&\left(a +b\lim _{n \to \infty}\int_{\Omega}|\nabla u_{n}|^p dx\right)^{s\cdot{\rm sgn}b}\nonumber\\
   & & \quad \left[\lim _{n \to \infty}\int_{\Omega}u_{n}|\nabla u_{n} |^{p-2}\nabla u_{n} \cdot\nabla\phi dx + \lim _{n \to \infty}\int_{\Omega}\phi|\nabla u_{n}|^p dx\right]
\end{eqnarray}
Now, by H\"{o}lder's inequality and weak convergence,
\begin{eqnarray*}
  0 &\leq&  \left(a +b\lim _{n \to \infty}\int_{\Omega}|\nabla u_{n}|^p dx\right)^{s\cdot{\rm sgn}b} \lim _{n \to \infty}\left|\int_{\Omega}u_{n}|\nabla u_{n} |^{p-2}\nabla u_{n} \cdot\nabla\phi dx\right|\\
     &\lesssim& \lim _{n \to \infty}\left|\int_{\Omega}u_{n}|\nabla u_{n} |^{p-2}\nabla u_{n} \cdot\nabla\phi dx\right|\\
      &\leq& \lim _{n \to \infty}\left(\int_{B\left(x_{k}, 2 \varepsilon\right) \cap \Omega}|\nabla u_{n}|^p dx\right)^{\frac{1}{p'}}    \lim _{n \to \infty}\left(\int_{B\left(x_{k}, 2 \varepsilon\right) \cap \Omega}| u_{n}|^p |\nabla\phi|^p dx\right)^{\frac{1}{p}}\\
      &\lesssim&\left(\int_{B\left(x_{k}, 2 \varepsilon\right) \cap \Omega}| u|^{p}|\nabla \phi|^{p} d x\right)^{\frac{1}{p}}\\
      &\leq&\left(\int_{B\left(x_{k}, 2 \varepsilon\right) \cap \Omega}|\nabla \phi|^{N} d x\right)^{\frac{p}{N}}\left(\int_{B\left(x_{k}, 2 \varepsilon\right) \cap \Omega}| u|^{p^*} d x\right)^{\frac{p}{p^*}}\\
      &\lesssim&\left(\int_{B\left(x_{k}, 2 \varepsilon\right) \cap \Omega}| u|^{p^*} d x\right)^{\frac{p}{p^*}}\to 0, \text{ as } \varepsilon \to 0.
\end{eqnarray*}
Substituting the above inequality into \eqref{Eq4.1}, we have
\begin{eqnarray*}
  0&=&\lim_{\varepsilon\to 0}\left[ - \int_{\Omega} \phi d \eta +  \int_{\Omega} g(x, u) u \phi d x - \lim _{n \to \infty} \left(a +b\int_{\Omega}|\nabla u_{n}|^p dx\right)^{s\cdot{\rm sgn}b}\int_{\Omega}|\nabla u_{n}|^{p-2}\nabla u_{n} \cdot\nabla(u_{n} \phi) dx\right] \\
   &=& - \eta_{k} - \left(a +b\lim _{n \to \infty}\int_{\Omega}|\nabla u_{n}|^p dx\right)^{s\cdot{\rm sgn}b}\lim_{\varepsilon\to 0}\left[ \lim _{n \to \infty}\int_{\Omega}\phi|\nabla u_{n}|^p dx\right]\\
   &\leq& - \eta_{k} - a^{s\cdot{\rm sgn}b}\mu_{k},
\end{eqnarray*}
which implies $\eta_{k}+ a^{s\cdot{\rm sgn}b}\mu_{k} \leq 0$. Therefore, $\eta_{k}=0$ and $\mu_{k}=0$ for any $k$. Consequently, $\|u_n\|\to \|u\|$.

\end{proof}

\begin{lemma}\label{Lem4.3}
Assume that conditions (\~{f}$_1$) is sufficed. Then there exists $\rho>0$ such that for all $u \in W_0$ with $\|u\| \leq \rho$ we have $I_{b}(u) \geq 0$ and $I_{b}(e) \leq 0$ for all $e \in \mathbb{R}$ with $|e| \leq \rho$.
\end{lemma}
\begin{remark}
(\~{f}$_1$) implies (f$_2$).  Since the proof of Lemma \ref{Lem3.3} uses only the continuity of the embedding and does not use the tightness of the embedding, the conclusion is still valid when the growth condition is promoted from the subcritical case to the critical case. And the proof method is similar, so we omitting the proof of Lemma \ref{Lem4.3}.
\end{remark}

\begin{proof}[Proof of Theorem \ref{Thm2}]
Combining  Lemma \ref{Lem2.1} with   Lemma \ref{Lem4.1},  Lemma \ref{Lem3.2} and  Lemma \ref{Lem4.3}, we obtain the proof of Theorem \ref{Thm2}.
\end{proof}

\section*{Acknowledgments}
The author wishes to express his thanks to Professor Peihao Zhao from the School of
Mathematics and Statistics in Lanzhou University for giving a guide to nonlinear functional
analysis.  




\end{document}